\documentclass[journal,9pt]{IEEEtran}

\usepackage{pifont}
\usepackage{bbding}  

\IEEEoverridecommandlockouts                  

\usepackage{float}
\usepackage{amsmath,amssymb,epsf,epsfig,times}
\usepackage[all]{xy}
\usepackage{graphicx,color}
\usepackage{subfigure}
\usepackage{url}
\usepackage{cite}
\usepackage{epstopdf}
\usepackage{epsfig}
\usepackage{mathrsfs} 

\newtheorem{theorem}{Theorem}[section]
\newtheorem{lemma}{Lemma}[section]

\def\QED{\mbox{\rule[0pt]{1.5ex}{1.5ex}}}
\def\endproof{\hspace*{\fill}~\QED\par\endtrivlist\unskip}
\newcommand{\re}{\mathbb{R}}

\newtheorem{definition}[theorem]{Definition}
\newtheorem{assumption}[theorem]{Assumption}

\newtheorem{remark}[theorem]{Remark}
\newtheorem{corollary}[theorem]{Corollary}

\newcommand{\OMIT}[1]{}
\allowdisplaybreaks[4]

		\title{Distributed Average Tracking of Heterogeneous Physical Second-order Agents With No Input Signals Constraint}
		
		\author{Sheida Ghapani, Student Member, IEEE, Salar Rahili, Student Member, IEEE, Wei Ren, Fellow, IEEE
		\thanks{Sheida Ghapani, Salar Rahili and Wei Ren are with the Department of Electrical and Computer Engineering, University of California, Riverside, CA, 92521, USA.
			{Email: sghap001@ucr.edu, srahi001@ucr.edu, ren@ee.ucr.edu.} 
			}
		}
\begin{document}
		\maketitle	
		
\begin{abstract}		
This paper addresses distributed average tracking of physical second-order agents with heterogeneous nonlinear dynamics, where there is no constraint on input signals.
The nonlinear terms in agents' dynamics are heterogeneous, satisfying a Lipschitz-like condition that will be defined later and is more general than the Lipschitz condition.
In the proposed algorithm, a control input and a filter are designed for each agent.
Each agent's filter has two outputs and the idea is that the first output  estimates the average of the input signals and the second output estimates the average of the input velocities asymptotically.
In parallel, each agent's position and velocity are driven to track, respectively, the first and the second outputs. 
Having heterogeneous nonlinear terms in agents' dynamics necessitates designing the filters for agents. Since the nonlinear terms in agents' dynamics can be unbounded and the input signals are arbitrary, novel state-dependent time-varying gains are employed in agents' filters and control inputs to overcome these unboundedness effects.
Finally the results are improved to achieve the distributed average tracking for a group of double-integrator agents, where there is no constraint on input signals and the filter is not required anymore.
Numerical simulations are also presented to illustrate the theoretical results.	
		\end{abstract}
\IEEEpeerreviewmaketitle	
	\section{Introduction}
In this paper,  the  distributed average tracking problem for a team of agents is studied, where each agent uses local information and local interaction to calculate the average of individual time-varying input signals, one per agent.
The problem has found applications in distributed sensor fusion \cite{olfati2004consensus}, feature-based map merging \cite{aragues2012distributed}, and distributed Kalman filtering \cite{bai2011}, where the scheme has been mainly used as an estimator.
However, there are some applications such as region following formation control \cite{cheah2009region} or distributed continuous-time convex optimization  \cite{SalarTac} that require the agents' physical
states instead of estimator
states to converge to a time-varying network quantity, where each agent only has a local and incomplete copy of that quantity.
Since the desired trajectory (the average of individual input signals) is time varying and not available to any agent, distributed average tracking poses more theoretical challenges than the consensus and distributed tracking problems.

In the literature, some researchers have employed linear distributed algorithms  
for some groups of time-varying input signals that are required to satisfy restrictive constraints \cite{spanos2005dynamic,freeman2006,bai2010,Zhumartinez2010,montijano2014robust,kiaDACsingularity}.
In \cite{freeman2006}, a proportional algorithm and a proportional-integral algorithm are proposed to achieve distributed average tracking for slowly-varying
input signals with a bounded tracking error.
In \cite{bai2010}, an extension of the proportional integral algorithm is employed for a special group of time-varying input signals with a common denominator in their Laplace transforms, where the denominator is required to be used in the estimator  design.
In  \cite{Zhumartinez2010} and \cite{montijano2014robust}, a discrete-time distributed average tracking problem is addressed, where \cite{montijano2014robust}  extends the proposed algorithm in \cite{Zhumartinez2010} by introducing  a time-varying sequence of damping factors to achieve robustness to initialization errors.

Note that the linear algorithms cannot ensure distributed average tracking for a more general group of input signals. Therefore, some researchers employ nonlinear tracking algorithms.
In \cite{Nosrati2012}, a class of nonlinear algorithms is proposed for input signals with bounded deviations, where the tracking error is proved to be bounded.
A non-smooth algorithm is proposed in \cite{chen2012distributed}, which is able to track time-varying input signals with bounded derivatives.
In the aforementioned works, the distributed average tracking problem is studied from a distributed estimation perspective without the requirement for agents to obey certain physical dynamics.
However, there are various applications, where the distributed average tracking problem is relevant for designing distributed control laws for physical agents.
The region-following formation control is one application \cite{cheah2009region}, where a swarm of robots move inside a dynamic region while keeping a desired formation.
In these applications, the dynamics of the physical agents must be taken into account in the control law design, where the dynamics themselves introduce further challenges in tracking and stability analysis.
Thus, distributed average tracking for physical agents with linear dynamics is studied in \cite{fei2015DAC,zhao2013distributed,FeiRobust,ghapani2015distributed}.
Refs.  \cite{fei2015DAC} studies the distributed average tracking for 
  input signals with bounded accelerations.
 Ref. \cite{ghapani2015distributed} introduces a discontinuous algorithm and filter for a group of physical double-integrator agents, where each agent uses the relative positions and neighbors' filter outputs to remove the velocity measurements.
 

However, in real applications physical agents might need to track the average of a group of time-varying input signals, where both physical agents and input signals have more complicated dynamics than single- or double-integrator dynamics.
Therefore, the control law designed for physical agents with linear dynamics can no longer be used directly for physical agents subject to more complicated dynamic equations.	
For example, \cite{fei2015EulerDAC} extends a proportional-integral control scheme  to achieve distributed average tracking for physical Euler-Lagrange systems for two different kinds of input signals with steady states and with bounded derivatives.

In this paper, an algorithm is introduced to achieve the distributed
average tracking for physical second-order agents with heterogeneous
nonlinear dynamics, where there is no constraint on the input signals.
Here, the nonlinear terms in agents' dynamics are heterogeneous,
satisfying a Lipschitz-like condition that will be defined later and
is more general than the Lipschitz condition. Therefore, the agents'
dynamics can cover many well-known systems such as car-like
robots.Due to the presence of the nonlinear heterogeneous terms
in the agents' dynamics, a local filter is introduced for each agent to
estimate the average of input signals and input velocities. Note that
the unknown terms in agents' dynamics can be unbounded and the
input signals are arbitrary, which create extra challenges. Therefore,
a time-varying state-dependent gains are employed in the filter's
dynamics and control law to overcome the unboundedness effects.In
the special case, where the agents have double-integrator dynamics,
the filter is not required anymore. Thus, the algorithm is modified to
drive the agents' positions and velocities to directly track the average
of the input signals and the input velocities. Here, by introducing
novel time-varying state-dependent gains in the control law, the
distributed average tracking is achieved for a group of arbitrary input
signals. A preliminary version of the work has appeared in \cite{ghapani2016distributed}. The
current work contains more rigorous algorithm design and proofs,
an additional result on heterogeneous double-integrator agents, and
numerical results, which is consistent with the TAC submission
policy.

The advantages of the proposed algorithms in comparison with the
literature are summarized as follows.
\begin{enumerate}
	\item In the first part of this paper, the agents are described by physical  heterogeneous nonlinear second-order dynamics.
	 The nonlinear terms in agents' dynamics, satisfying a Lipschitz-like condition that will be defined later, are heterogeneous and hence more general and more realistic. 
	In contrast, in \cite{chen2012distributed,fei2015DAC}, the agents' dynamics are assumed to be homogeneous single- or double-integrator dynamics. 
	The results for single- or double-integrator dynamics are not applicable to physical heterogeneous nonlinear  agents.
	
	\item  
	In the existing algorithms in the literature, there are always constraints on input signals.
	For example, in \cite{fei2015DAC}, the second derivative of the input signals are assumed to be bounded.
	In contrast, by proposing a new distributed control law in combination with  a local filter for each agent,  there is no limitation on the input signals in the current paper.
	The novelty of the local filters is that by introducing new time-varying state-dependent gains, 
	the distributed average tracking problem can be achieved for an arbitrary group of input signals. 
\end{enumerate}
	\section{Notations and Preliminaries}	
	\subsection{Notations}
	The following notations are adopted throughout this paper.
	$\mathbb{R}$ denotes the set of all real numbers and $\mathbb{R}^+$ denotes the set of all positive real numbers. 
	The transpose of matrix $A$ and vector $x$ are shown as $A^T$ and $x^T$, respectively. 
	Let $\mathbf{1}_n$ and $\mathbf{0}_n$ denote, respectively, the $n \times 1$ column vector of all ones and all zeros.
	Let $\mbox{diag}(z_1,\ldots,z_p)$ be the diagonal matrix with diagonal entries $z_1$ to $z_p$.
	We use 
	$\otimes$ to denote the Kronecker product, and $\mbox{sgn}(\cdot)$ to denote the $\mbox{signum}$ function defined component-wise. 
	For a vector function ${x(t):\re\mapsto\re^m}$, define $\|x\|_\mathfrak{p}$ as the $\mathfrak{p}$-norm; 
	${x(t)\in\mathbb{L}_2}$ if
	$\int_{0}^{\infty} \|x(\tau)\|_2^2 \mbox{d}\tau<\infty$ and
	${x(t)\in\mathbb{L}_{\infty}}$ if for each element of $x(t)$, denoted as
	$x_i(t)$, ${\sup_{t \geq 0}|x_i(t)|<\infty}$, $i=1,\ldots,m$.
	
	\subsection{Graph Theory}
	An \textit{undirected} graph $G \triangleq (V,E)$ is used to characterize the interaction topology among the agents, where ${V \triangleq \{1,\ldots,n\}}$ is the node set and $E \subseteq V \times V$ is the edge set.
	An edge $(j,i) \in E$ means that node $i$ and node $j$ can obtain information from each
	other and they are neighbors of each other.
	Self edges $(i,i)$ are not considered here.
	The set of neighbors of node $i$ is denoted as $N_i$.
	The \textit{adjacency matrix} ${\mathbf{A} \triangleq [a_{ij}] \in \mathbb{R}^{n \times n}}$ of the graph $G$ is defined such that the edge weight ${a_{ij}=1}$ if ${(j,i) \in E}$ and ${a_{ij}=0}$ otherwise. For an undirected graph, ${a_{ij}=a_{ji}}$.
	The \textit{Laplacian matrix} ${L \triangleq [l_{ij}] \in \mathbb{R}^{n \times n}}$ associated with $\mathbf{A}$ is defined as ${l_{ii}=\sum_{j \ne i} a_{ij}}$ and ${l_{ij}=-a_{ij}}$, for ${i \ne j}$.
	For an undirected graph, $L$ is symmetric positive semi-definite.
	Let $m$ denote the number of edges in $E$, where the edges $(j,i)$ and $(i,j)$ are counted only once. By arbitrarily assigning an orientation for the edges in $G$, let $D \triangleq [d_{ij}] \in  \mathbb{R}^{n \times m}$ be the \textit{incidence matrix} associated with $G$, where $d_{ij} = -1$ if the edge $e_j$ leaves node $i$, $d_{ij} = 1$ if it enters node $i$, and $d_{ij} = 0$ otherwise.
	The \textit{Laplacian matrix} $L$ is then given by $L=DD^T$ \cite{GodsilRoyle01}.
	
	\subsection{Nonsmooth Analysis}\vspace{-0.4cm}
	Consider a vector-valued differential equation
	\begin{align} \label{dif-equ}
		\dot{x}(t)=f(t,x(t)),
	\end{align}
	where $	t \in \mathbb{R}$ and $x(t) \in \mathbb{R}^n$.
	\begin{definition}
		
		For the differential equation \eqref{dif-equ}, define the Filippov set-valued map $\mathcal{K} [f](t,x(t)) \triangleq \cap_{\delta >0} \cap_{u(N)=0} \bar{co} \Big (f(t, B(x(t),\delta)-N) \Big)$, where $\bar{co}$ denotes the convex closure, $\cap_{u(N)=0}$ denotes the intersection over all
		sets of Lebesgue measure zeros and $B(x,\delta)$ is the open ball of radius $\delta$ centered at $x$ \cite{filippov1960differential}.
	\end{definition}
	\begin{definition}
		Replace the differential equation \eqref{dif-equ} with the differential inclusion
		\begin{align}\label{inc}
			\dot{x}(t) \in \mathcal{K} [f](t,x(t)).
		\end{align}
		A vector function $x(\cdot): \mathbb{R} \to \mathbb{R}^n$ is called a Filippov solution of \eqref{dif-equ} on $[t_0, t_1]$, $t_1 \leq \infty$, if $x(\cdot)$ is absolutely continuous and satisfies \eqref{inc} for almost all $t \in [t_0, t_1]$.
	\end{definition}
	\begin{lemma}\cite{filippov1960differential}\label{fil-exs}
		Suppose that $f(t,x(t))$ in \eqref{dif-equ} is measurable and locally essentially bounded,  that is, bounded on a bounded neighborhood of every point excluding sets of measure zero. Then, for any $x(0) \in \mathbb{R}^n$, there exists a Filippov solution of \eqref{dif-equ} with initial condition $x(0)=x_0$. 
	\end{lemma}
	Let $\mathcal{W}[(x(t))] : \mathbb{R}^n \to \mathbb{R}$ be a locally Lipschitz function of $x(t)$.
	The generalized gradient of $\mathcal{W}[x(t)]$ is defined 
	\begin{align*}
		\partial \tilde{\mathcal{W}} \triangleq co \Big \{ \lim\limits_{w \to x} \nabla \mathcal{W}(w) : w \in  \mathbb{R}^n, w \notin \Omega_{\mathcal{W}} \cup \mathcal{S}  \Big \},
	\end{align*}
	where $co(\cdot)$ is the convex hull, $\Omega_{\mathcal{W}}$ is the set of points in which $\mathcal{W}[x(t)]$ is not differentiable
	and $\mathcal{S}$ is a set of measure zero that can be arbitrarily chosen so as to simplify the calculation. 
	The set-valued Lie derivative of
	$\mathcal{W}[x(t)]$ with respect to $x(t)$, the trajectory of \eqref{dif-equ}, is defined as $\dot{\tilde{\mathcal{W}}} \triangleq \cap_{\zeta  \in \partial \tilde{\mathcal{W}} } \zeta^T \mathcal{K} [f]$.
	\begin{assumption} \label{conn-graph}
		Graph $G$ is connected.
	\end{assumption}
	\begin{lemma} \cite{GodsilRoyle01} \label{eigen}
		Under Assumption \ref{conn-graph}, the \textit{Laplacian matrix} $L$ has a simple zero eigenvalue such that $0=\lambda_1(L)<\lambda_2(L) \leq \ldots \leq \lambda_n(L)$, where $\lambda_i(\cdot)$ denotes the $i$th eigenvalue. Furthermore, for any vector $y \in \mathbb{R}^n$ satisfying ${\mathbf{1}_n^T y=0}$, we have $\lambda_2(L) y^Ty \leq y^T L y \leq \lambda_n(L) y^Ty$.
	\end{lemma}
		\begin{lemma}\label{nonsgn}
		For any vector $x \in \mathbb{R}^n$, we have
			\begin{align}
x^T L DW \mbox{sgn} (D^T x) \geq &	\lambda_2 (L) x^T D W \mbox{sgn}(D^Tx), \label{wght-D}		
			\end{align}
				\end{lemma}
where $W$ is a positive-definite diagonal matrix.
 
\emph{Proof}:
If $D^T x =\mathbf{0}_n $, \eqref{wght-D} holds. 
However, if $D^T x \neq \mathbf{0}_n$, then we replace $L$ with $D D^T$ on the left side of \eqref{wght-D}. Thus, we will have
\begin{align*}
 x^T L DW \mbox{sgn} (D^T x)  &=
 x^T D D^T DW \mbox{sgn} (D^T x) \\ &=(D^T x)^T D^T D W \mbox{sgn} (D^T x).
\end{align*}
Note that both $D^T  D$ and $D D^T$ have the same set of nonzero eigenvalues $\Lambda=\{\lambda_m,\lambda_{m+1}, \ldots, \lambda_n \}$ \cite{zelazo2007agreement}.
Suppose that $S$ is the space spanned by the eigenvectors belonging to the nonzero eigenvalues of $D^T  D$.
If $D^T x \in S$, there is a $\lambda_i(D^TD) \in \Lambda$ such that $ D^TD (D^T x)= \lambda_i(D^TD) (D^T x)$.
Thus, we get that
\begin{align*}
 (D^T x)^T  D^T D W \mbox{sgn} (D^T x) \\ 
 =&  \lambda_i(D^TD) (D^T x)^T   W \mbox{sgn} (D^T x) \\
\geq & \lambda_2(D^T D)  (D^T x)^T   W \mbox{sgn} (D^T x) \\
=& \lambda_2(D^T D)  x^T D  W \mbox{sgn} (D^T x).
\end{align*}
If $D^T x \notin S$, it follows that $D^T x$ belongs to the null space of $D^TD$. 
Based on Lemma 3 in \cite{feidoubleintegrator}, the null space of the incidence matrix $D$ coincides with null space of $D^TD$.
Thus, $D (D^Tx)=\mathbf{0}_n$ which means $Lx=\mathbf{0}_n$.
It follows that $x$ belongs to the space spanned by vector $\mathbf{1}_n$ and hence $D^Tx=\mathbf{0}_n$.
This contradicts with $D^T x \neq \mathbf{0}_n$ and hence $D^T x \in S$.
%
\endproof
\section{Problem Statement}
Consider a multi-agent system consisting of $n$ physical agents described by the following heterogeneous nonlinear second-order dynamics	
\begin{align} \label{non-agents-dynamic}
	\dot{x}_i(t)=&v_i(t), \notag  \\
	\dot{v}_i(t)=&f_i(x_i(t),v_i(t),t)+u_i(t),  \qquad i=1,\ldots,n,
\end{align}
where $x_i(t)$, $v_i(t) $, $u_i(t) \in \mathbb{R}^\mathfrak{p}$ are the $i$th agent's position, velocity and control input, respectively, $f_i:\mathbb{R}^\mathfrak{p} \times \mathbb{R}^\mathfrak{p} \times \mathbb{R}^+ \to \mathbb{R}^\mathfrak{p}$ is a 
vector-valued nonlinear function which will be defined later.
Suppose that each agent has a time-varying input signal $x^r_i(t) \in \mathbb{R}^\mathfrak{p}$, $i=1,\ldots,n$, satisfying
\begin{align} \label{ref-dynamic}
	\dot{x}^r_i(t)= v_i^r(t),   \qquad
	\dot{v}_i^r(t)= a^r_i(t),
\end{align}
where $v_i^r(t)$, $a^r_i (t) \in \mathbb{R}^\mathfrak{p}$ are, respectively, the input velocity and the input acceleration.
Note that \eqref{ref-dynamic} is just used to show the relation between the input signal, input velocity and input acceleration. 

The goal is to design $u_i(t)$ for agent $i$, $i=1,\ldots,n$, to track the average of the input signals and input velocities, i.e.,
\begin{align}\label{goal}
	\lim \limits_{t \to \infty} \|x_i(t)-\frac{1}{n} \sum_{j=1}^n r_j(t)\|_2=& \mathbf{0}_{\mathfrak{p}}, \notag \\
	\lim \limits_{t \to \infty} \|v_i(t)-\frac{1}{n} \sum_{j=1}^n v_j^r(t)\|_2=& \mathbf{0}_{\mathfrak{p}}, 
\end{align}
where each agent has only local interaction with its neighbors and has access to only its own input signal, velocity, and acceleration. 	
As it was mentioned, there are applications, where the physical agents track the average of a group of time-varying signals while the physical agents  might be described by more complicated dynamics rather than linear dynamics.
Here, we investigate the distributed average tracking problem for a more general group of agents while there is no constraint on input signals.

%
\subsection{Distributed Average Tracking for Physical Heterogeneous Nonlinear Second-order Agents}\label{section2}
In this subsection, we study the distributed average tracking problem for a group of heterogeneous nonlinear second-order agents, where the nonlinear term $f_i(\cdot,\cdot,t)$ satisfies the Lipschitz-like condition and there is no constraint on input signals. 

\begin{assumption}\label{lip} 
	The vector-valued function $f_i(\cdot,\cdot,t)$ is continuous in $t$ and satisfies the following Lipschitz-like condition $\forall t \geq 0$
	\begin{align*}
	\left\{
	\begin{array}{rl}
		\|f_i(x,v,t)-f_i(y,z,t) \|_1 \leq & \rho_1 \|x-y \|_1  +\rho_2 \|v-z \|_1  +\rho_3,  \\
		\| f_i(\mathbf{0}_{\mathfrak{p}},\mathbf{0}_{\mathfrak{p}},t) \|_1 \leq & \rho_4,
		\end{array} \right.
	\end{align*}
	where $x$, $v$, $y$, $z \in \mathbb{R}^\mathfrak{p}$, and $\rho_1$, $\rho_2$, $\rho_3$, $\rho_4 \in \mathbb{R}^+$.
\end{assumption}
	\begin{remark}
	Note that Assumption \ref{lip} is more general than the Lipschitz condition, satisfied by many well-known systems such as the pendulum system with a control torque, car-like robots, the Chua's circuit, the Lorenz system, and the Chen system \cite{mei2013distributed}.
	In fact, the term $f_i(\cdot,\cdot,t)$ is general enough to represent both the nonlinear dynamics and possible bounded disturbances.
\end{remark}

It should be noted that the nonlinear term $f_i(\cdot,\cdot,t)$ is unknown and the input acceleration $a^r_i(t)$ is arbitrary and can be unbounded.
Therefore, a novel filter is introduced with time-varying state-dependent  gains to estimate the average of input signals and input velocities in the presence of these challenges. For notational simplicity, we will remove the index $t$ from the variables in the reminder of the paper.
Consider the following local filter for agent $i$  
\begin{align}
		p_i=& z_i+x^r_i, \notag \\
		\ddot{z}_i=& -\kappa \big ( p_i-x^r_i \big) -\kappa \big ( q_i-v^r_i \big)  \notag \\
		&-  \beta\sum\limits_{j=1}^{n} a_{ij}   (\psi_{i}+\psi_{j}) \mbox{sgn} \big[
		(p_i+q_i)-(p_j+q_j) \big], \label{p-term-adap}
\end{align}
where $p_i$, $q_i \in \mathbb{R}^\mathfrak{p}$ are the filter outputs, $q_i=\dot{p}_i$, $z_i \in \mathbb{R}^\mathfrak{p}$ is an auxiliary filter variable, $\psi_{i}=\|x^r_i\|_1+\|v^r_i\|_1+\|a^r_i\|_1+\gamma$, and $\kappa, \beta$, $\gamma \in \mathbb{R}^+$ will be designed later.

The control input $u_i$ is designed as
\begin{align}\label{con-inp-adap}
	u_i=&-\eta \psi'_i \big[\tilde{x}_i + \tilde{v}_i + \mbox{sgn}(\tilde{x}_i+\tilde{v}_i) \big]  +\dot{q}_i,  
\end{align}
where $\tilde{x}_i=x_i-p_i$, $\tilde{v}_i=v_i-q_i$, $\psi'_i=\|x_i\|_1+\|v_i\|_1+\gamma$ and $\eta \in \mathbb{R}^+$ to be designed.
Here, the time-varying state-dependent  gains $\psi_i$ and $\psi'_i$ are employed in the filter's dynamics \eqref{p-term-adap} and the control law \eqref{con-inp-adap} to overcome the unboundedness challenges of the heterogeneous unknown term $f_i(\cdot,\cdot,t)$ and the arbitrary input acceleration $a^r_i$.
\begin{theorem} \label{DAC-lip-adap-adap}
	Under the control algorithm given by  \eqref{p-term-adap}-\eqref{con-inp-adap} for the system \eqref{non-agents-dynamic}, the distributed average tracking goal \eqref{goal} is achieved asymptotically, provided that Assumptions \ref{conn-graph} and \ref{lip} hold and the control gains satisfy the constraints $\beta>\kappa$, $\kappa >\max \{1, \frac{\lambda^2_{\max} (L)}{2 \lambda_2^2(L)}\}$, $\gamma > \rho_3+\rho_4$, and $\eta> \max \{ 1, \rho_1,\rho_2 \}$.
\end{theorem}
\emph{Proof}:
First, it is proved that for $i=1,\cdots,n,$ 
\begin{align} \label{pi-avrri}
\lim\limits_{t \to \infty} p_i = \frac{1}{n} \sum\limits_{j=1}^{n} r_j, \quad \lim\limits_{t \to \infty} q_i = \frac{1}{n} \sum\limits_{j=1}^{n} v^r_j. \end{align} 	
Using $q_i=\dot{p}_i$, the local filter's dynamics \eqref{p-term-adap} can be rewritten as
\begin{align}\label{fil-clos}
	\dot{p}_i=& q_i,  \\
	\dot{q}_i=& -\kappa \big ( p_i-x^r_i \big) -\kappa \big ( q_i-v^r_i \big)  \notag \\
	&-  \beta\sum\limits_{j=1}^{n} a_{ij}   (\psi_{i}+\psi_{j}) \mbox{sgn} \big[
	(p_i+q_i)-(p_j+q_j) \big]+a^r_i. \notag
\end{align}
Due to the existence of the signum function in the algorithm \eqref{p-term-adap}, the closed-loop dynamics \eqref{fil-clos} is discontinuous.
Therefore, the solutions should be investigated in terms of differential inclusions by using nonsmooth analysis \cite{filipi,cortes2008discontinuous}.
Since the signum function is measurable and locally essentially bounded, the Filippov solutions for the closed-loop
dynamics \eqref{fil-clos} always exist and are absolutely continuous.
Let $r=[r^T_1,\ldots,r^T_n ]^T$, $v^r=[{v_1^r}^T,\ldots,{v_n^r}^T ]^T$, $a^r=[{a_1^r}^T,\ldots,{a_n^r}^T ]^T$, $p=[p^T_1,\ldots,p^T_n ]^T$, $q=[q^T_1,\ldots,q^T_n ]^T$, $\tilde{p}_i=p_i-\frac{1}{n}\sum\limits_{j=1}^{n} p_j$ and $\tilde{q}_i=q_i-\frac{1}{n}\sum\limits_{j=1}^{n} q_j$.
Defining
$M=I_n- \frac{1}{n}\textbf{1}^T_n \textbf{1}_n$, $\tilde{p}=(M \otimes I_\mathfrak{p}) p$ and $\tilde{q}=(M \otimes I_\mathfrak{p}) q$, we get that
\begin{align}\label{cls-inc}
	\begin{bmatrix}
		\dot{\tilde{p}} \\
		 \dot{\tilde{q}}
	\end{bmatrix} \in^{a.e.} \mathcal{K} [\mathfrak{f}](\tilde{p},\tilde{q})
\end{align}
where $a.e.$ stands for ``almost everywhere", 
${\mathfrak{f}=
	\begin{bmatrix}
	\mathfrak{f}_p^T & \mathfrak{f}_q^T 
	\end{bmatrix}}^T$ and
\begin{align*}
\mathfrak{f}_{\tilde{p}}=& \tilde{q}, \notag \\
\mathfrak{f}_{\tilde{q}}=& -\kappa \tilde{p} + \kappa (M \otimes I_\mathfrak{p})r-\kappa \tilde{q} + \kappa (M \otimes I_\mathfrak{p}) v^r \\
&- \beta(DW  \otimes I_\mathfrak{p})  \mbox{sgn} [( D^T \otimes I_\mathfrak{p})(\tilde{p}+\tilde{q})] + (M  \otimes I_\mathfrak{p}) a^r,
\end{align*}
where $W$ is a diagonal matrix and $MD=D$.
The $k$th diagonal element of the matrix $W$ describes the $k$th edge weight. If the $k$th edge is between node $i$ and node $j$, then it is equal to $\psi_i+\psi_j$.

Consider the following Lyapunov function candidate
\begin{align} \label{V1}
	V_1=& \frac{1}{2} 
	\begin{bmatrix}
		\tilde{p}^T && \tilde{q}^T
	\end{bmatrix}
	\big ( L \otimes
		\begin{bmatrix}
		2\kappa  I_{\mathfrak{p}}  &&  I_{\mathfrak{p}}  \\
		 I_{\mathfrak{p}}   && I_{\mathfrak{p}}
	\end{bmatrix}
	\big )
	\begin{bmatrix}
		\tilde{p} \\
		\tilde{q}
	\end{bmatrix}.
\end{align}
Since $(\mathbf{1}_n \otimes I_\mathfrak{p}) ^T \tilde{p} =\mathbf{0}_{n\mathfrak{p}}$ and $(\mathbf{1}_n \otimes I_\mathfrak{p}) ^T \tilde{q} =\mathbf{0}_{n\mathfrak{p}}$, by using Lemma \ref{eigen}, we will have
\begin{align*}
V_1 \geq &\frac{1}{2}  \begin{bmatrix}
\tilde{p}^T && \tilde{q}^T 
\end{bmatrix} 
\begin{bmatrix}
2\kappa \lambda_2(L)  I_{n\mathfrak{p}}  && L \otimes I_{n\mathfrak{p}} \\
L \otimes I_{n\mathfrak{p}}  && \lambda_2(L) I_{n\mathfrak{p}}
\end{bmatrix}
\begin{bmatrix}
\tilde{p} \\
\tilde{q}
\end{bmatrix}.
\end{align*}
 If $\kappa >\frac{\lambda^2_{\max}(L)}{2\lambda^2_2(L)} $, then it can be proved that the matrix $\begin{bmatrix}
2\kappa \lambda_2(L)  I_{n\mathfrak{p}}  && L \otimes I_{n\mathfrak{p}} \\
L \otimes I_{n\mathfrak{p}}  && \lambda_2(L) I_{n\mathfrak{p}}
\end{bmatrix}$ and hence $V_1$ are positive definite. 
Since $V_1$ is a continuous function, its set-valued Lie derivative along \eqref{cls-inc} is given as
\begin{align}\label{lypV1}
	\dot{\tilde{V}}_1=& \mathcal{K} \Big[ 2 \kappa  \tilde{p}^T (L \otimes I_\mathfrak{p}) \tilde{q}+\tilde{q}^T (L \otimes I_\mathfrak{p}) \tilde{q} - \kappa \tilde{p}^T(L  \otimes I_\mathfrak{p})\tilde{p}\notag \\
	& +\kappa \tilde{p}^T(L  \otimes I_\mathfrak{p}) r - \kappa \tilde{p}^T(L  \otimes I_\mathfrak{p})\tilde{q} + \kappa \tilde{p}^T(L  \otimes I_\mathfrak{p}) v^r\notag \\
	&- \beta \tilde{p}^T ( L D W  \otimes I_\mathfrak{p})  \mbox{sgn} [( D^T \otimes I_\mathfrak{p})(\tilde{p}+\tilde{q})]\notag \\
	&+\tilde{p}^T (L \otimes I_\mathfrak{p})  a^r - \kappa \tilde{q}^T(L  \otimes I_\mathfrak{p})\tilde{p} +\kappa \tilde{q}^T(L  \otimes I_\mathfrak{p}) r \notag \\
	&- \kappa \tilde{q}^T(L  \otimes I_\mathfrak{p})\tilde{q} + \kappa \tilde{q}^T(L  \otimes I_\mathfrak{p}) v^r\notag \\
	&- \beta \tilde{q}^T ( L D W  \otimes I_\mathfrak{p})  \mbox{sgn} [( D^T \otimes I_\mathfrak{p})(\tilde{p}+\tilde{q})]\notag \\
	&+\tilde{q}^T (L \otimes I_\mathfrak{p})  a^r \Big ] \notag \\
	=&- \kappa \tilde{p}^T (L \otimes I_\mathfrak{p})  \tilde{p} - ( \kappa-1) \tilde{q}^T  ( L \otimes I_\mathfrak{p} )  \tilde{q} \notag \\
	&- \beta (\tilde{p}^T+\tilde{q}^T) ( L D W  \otimes I_\mathfrak{p})  \mathcal{K} \Big[ \mbox{sgn} [( D^T \otimes I_\mathfrak{p})(\tilde{p}+\tilde{q})] \Big ] \notag \\
	&+(\tilde{p}^T +\tilde{q}^T ) (L \otimes I_\mathfrak{p})  (\kappa r +\kappa v^r +a^r).
\end{align} 
If $( D^T \otimes I_\mathfrak{p})(\tilde{p}+\tilde{q})=\mathbf{0}_{m\mathfrak{p}}$, then $( L \otimes I_\mathfrak{p})(\tilde{p}+\tilde{q})=\mathbf{0}_{n\mathfrak{p}}$.
Therefore, it is concluded that $- \beta (\tilde{p}^T+\tilde{q}^T) ( L D W  \otimes I_\mathfrak{p})  \mathcal{K} \Big[ \mbox{sgn} [( D^T \otimes I_\mathfrak{p})(\tilde{p}+\tilde{q})] \Big ]=0$.
If $( D^T \otimes I_\mathfrak{p})(\tilde{p}+\tilde{q}) \neq  \mathbf{0}_{m\mathfrak{p}}$, then $\mathcal{K} \Big[ \mbox{sgn} [( D^T \otimes I_\mathfrak{p})(\tilde{p}+\tilde{q})] \Big ]=\Big \{ \mbox{sgn} [( D^T \otimes I_\mathfrak{p})(\tilde{p}+\tilde{q})] \Big \}$.
Thus, in both cases the set-valued Lie derivative of $V_1$ is a singleton. Note that the function $V_1$ is continuously differentiable.
By using Lemma \ref{nonsgn}, it follows from $\dot{V}_1 \in \dot{\tilde{V}}_1$, where $\dot{V}_1$ denotes the derivative of
$V_1$, that 
\begin{align*}
	\dot{V}_1
	\leq & - \kappa \tilde{p}^T (L \otimes I_\mathfrak{p})  \tilde{p} - ( \kappa-1) \tilde{q}^T  ( L \otimes I_\mathfrak{p} )  \tilde{q} \notag \\
	&- \beta \lambda_2(L) (\tilde{p}^T +\tilde{q}^T) (  D W  \otimes I_\mathfrak{p})  \mbox{sgn} [( D^T \otimes I_\mathfrak{p})(\tilde{p}+\tilde{q})] \notag \\
&+(\tilde{p}^T +\tilde{q}^T ) (L \otimes I_\mathfrak{p})  (\kappa r +\kappa v^r +a^r)\\
	= &- \kappa \tilde{p}^T (L \otimes I_\mathfrak{p})  \tilde{p} - ( \kappa-1) \tilde{q}^T  ( L \otimes I_\mathfrak{p} )  \tilde{q} \notag \\
	& - \beta \sum\limits_{i=1}^{n}  \sum\limits_{j=1}^{n} a_{ij} (\psi_i+\psi_j) \Big [ (\tilde{p}_i+\tilde{q}_i)-(\tilde{p}_j+\tilde{q}_j) \Big ]^T  \times  \notag \\
	&\mbox{sgn} \Big[ (\tilde{p}_i+\tilde{q}_i)-(\tilde{p}_j+\tilde{q}_j) \Big ]
	\notag \\
	&+(\tilde{p}^T +\tilde{q}^T ) (L \otimes I_\mathfrak{p})  (\kappa r +\kappa v^r +a^r)  \\
	=&- \kappa \tilde{p}^T (L \otimes I_\mathfrak{p})  \tilde{p} - ( \kappa-1) \tilde{q}^T  ( L \otimes I_\mathfrak{p} )  \tilde{q} \notag \\
	&-\beta \sum\limits_{i=1}^{n} \sum\limits_{j=1}^{n} a_{ij} (\psi_i+\psi_j)   \Big \| (\tilde{p}_i+\tilde{q}_i)-(\tilde{p}_j+\tilde{q}_j)  \Big \|_1 \\
	&+\sum\limits_{i=1}^{n} \Big [ \sum\limits_{j=1}^{n} a_{ij}  \Big \{(\tilde{p}_i+\tilde{q}_i)-(\tilde{p}_j+\tilde{q}_j) \Big \} \Big ]^T ( \kappa x^r_i+ \kappa v^r_i+a^r_i),
\end{align*}
where we have used \eqref{wght-D} to obtain the inequality and $p_i+q_i-p_j-q_j=\tilde{p}_i+\tilde{q}_i-\tilde{p}_j-\tilde{q}_j$ to obtain the second equality.
By using the triangular inequality, we can get that
\begin{align*}
	\dot{V}_1 \leq & - \kappa \tilde{p}^T (L \otimes I_\mathfrak{p})  \tilde{p} - ( \kappa-1) \tilde{q}^T  ( L \otimes I_\mathfrak{p} )  \tilde{q} \notag \\
	&-\beta \sum\limits_{i=1}^{n}  \sum\limits_{j=1}^{n} a_{ij} (\psi_i+\psi_j) \Big \| (\tilde{p}_i+\tilde{q}_i)-(\tilde{p}_j+\tilde{q}_j) \Big \|_1 \\
	&+\sum\limits_{i=1}^{n}  \sum\limits_{j=1}^{n} a_{ij} \Big\| (\tilde{p}_i+\tilde{q}_i)-(\tilde{p}_j+\tilde{q}_j)  \Big\|_1  \times \\
	& \big (\kappa \|x^r_i\|_1+ \kappa \|v^r_i\|_1+\|a^r_i\|_1 \big) \\
	=& - \kappa \tilde{p}^T (L \otimes I_\mathfrak{p})  \tilde{p} - ( \kappa-1) \tilde{q}^T  ( L \otimes I_\mathfrak{p} )  \tilde{q} \notag \\
	&+ \sum\limits_{i=1}^{n} (\kappa \|x^r_i\|_1+ \kappa \|v^r_i\|_1+\|a^r_i\|_1 -\beta\psi_i) \times \\
	&  \sum\limits_{j=1}^{n} a_{ij} \Big \|  (\tilde{p}_i+\tilde{q}_i)-(\tilde{p}_j+\tilde{q}_j)  \Big \|_1 \\
	&-\beta \sum\limits_{i=1}^{n}   \sum\limits_{j=1}^{n} a_{ij} \psi_j \Big \|  (\tilde{p}_i+\tilde{q}_i)-(\tilde{p}_j+\tilde{q}_j)  \Big \|_1.
	\end{align*}
	Since $\psi_{i}= \|x^r_i\|_1+  \|v^r_i\|_1+\|a^r_i\|_1+\gamma$ and $\beta > \kappa > 1 $, we will have 
	\begin{align}\label{dotV1}
	\dot{V}_1 \leq& - \kappa \tilde{p}^T (L \otimes I_\mathfrak{p})  \tilde{p} - ( \kappa-1) \tilde{q}^T  ( L \otimes I_\mathfrak{p} )  \tilde{q} \notag \\
	\leq & -\kappa \lambda_2(L) \tilde{p}^T  \tilde{p}-(\kappa -1) \lambda_2(L) \tilde{q}^T  \tilde{q} < 0,
\end{align}
where we have used $\kappa > 1$, Lemma \ref{eigen} and the fact that $(\mathbf{1}_n \otimes I_\mathfrak{p}) ^T \tilde{p} =\mathbf{0}_{n\mathfrak{p}}$ and $(\mathbf{1}_n \otimes I_\mathfrak{p}) ^T \tilde{q}=\mathbf{0}_{n\mathfrak{p}}$ to obtain the second inequality.
Using Theorem 4.10 in \cite{khalil2002nonlinear}, it is concluded that $\begin{bmatrix}
\tilde{p} \\
\tilde{q}
\end{bmatrix}=\mathbf{0}_{2n \mathfrak{p}}$ is globally exponentially stable, which means for $i=1,\ldots,n$, \begin{align}\label{pi-avqi}
\lim\limits_{t \to \infty} p_i =& \frac{1}{n} \sum\limits_{j=1}^{n} p_j, \quad \lim\limits_{t \to \infty} q_i = \frac{1}{n} \sum\limits_{j=1}^{n} q_j. 
\end{align}

Now, it is proved that $ \sum\limits_{j=1}^{n} p_j \to \sum\limits_{j=1}^{n} r_j$ and $ \sum\limits_{j=1}^{n} q_j \to \sum\limits_{j=1}^{n} v^r_j $.  
Defining the variables $S_p=\sum\limits_{j=1}^{n} (p_j - r_j)$ and
$S_q=\sum\limits_{j=1}^{n} (q_j - v^r_j)$,
we can get from \eqref{fil-clos} that
\begin{align*}
\begin{bmatrix}
\dot{S}_p \\
\dot{S}_q
\end{bmatrix} = \big (\begin{bmatrix}
0 && 1 \\
-\kappa  && - \kappa \\
\end{bmatrix} \otimes I_p \big)  \begin{bmatrix}
S_p \\
S_q
\end{bmatrix}= (A\otimes I_p \big) \begin{bmatrix}
S_p \\
S_q
\end{bmatrix}.
\end{align*}
If $\kappa>0$, the matrix $A$ is Hurwitz.
Therefore, $\lim\limits_{t \to \infty} \begin{bmatrix}
S_p \\
S_q
\end{bmatrix} =\mathbf{0}_{p} $, which means $\lim\limits_{t \to \infty} \sum\limits_{j=1}^{n} p_j = \sum\limits_{j=1}^{n} r_j$ and $\lim\limits_{t \to \infty} \sum\limits_{j=1}^{n} q_j = \sum\limits_{j=1}^{n} v^r_j$. 
Now, using \eqref{pi-avqi}, it is easy to see that \eqref{pi-avrri} holds. 

Second, it is proved that by using the control law \eqref{con-inp-adap} for \eqref{non-agents-dynamic}, $\lim\limits_{t \to \infty} x_i = p_i$ and $\lim\limits_{t \to \infty} v_i = q_i$ in parallel and hence it can be concluded that
$\lim\limits_{t \to \infty} x_i = \frac{1}{n} \sum\limits_{j=1}^{n} r_j$ and $\lim\limits_{t \to \infty} v_i = \frac{1}{n} \sum\limits_{j=1}^{n} v^r_j$. 
Define $\tilde{x}=[\tilde{x}_1^T,\ldots,\tilde{x}_n^T]^T$, $\tilde{v}=[\tilde{v}_1^T,\ldots,\tilde{v}_n^T]^T$, and $f(x,v,t)=[f_1^T(x_1,v_1,t),\ldots,f_n^T(x_n,v_n,t)]^T$.
Using the control law \eqref{con-inp-adap} for \eqref{non-agents-dynamic}, we get the closed-loop dynamics in vector form as
\begin{align}\label{cls-mn}
	\dot{\tilde{x}}=&\tilde{v}, \\
	\dot{\tilde{v}}=& f(x,v,t) -\eta (\tilde{x}+ \tilde{v}) -\eta ( \psi' \otimes I_{\mathfrak{p}}) \mbox{sgn}(\tilde{x}+\tilde{v}), \notag
\end{align}
 where $\psi' \triangleq \mbox{diag}(\psi'_{1} ,\ldots,\psi'_{n})$. Since the signum function is measurable and locally essentially bounded, the Filippov solution for the closed-loop
dynamics \eqref{cls-mn} exists.
Consider the Lyapunov function candidate
\begin{align}\label{LypV2}
V_2= \frac{1}{2} 
\begin{bmatrix}
\tilde{x}^T && \tilde{v}^T
\end{bmatrix}
\big( 
\begin{bmatrix}
2\eta  && 1 \\
1 && 1
\end{bmatrix}
\otimes I_{n\mathfrak{p}} \big)
\begin{bmatrix}
\tilde{x} \\
\tilde{v}
\end{bmatrix}.
\end{align}
It is easy to see that $V_2$ is positive definite.
By taking the set-valued Lie derivative of $V_2$, $\dot{\tilde{V}}_2$, along the Filippov set-valued map of \eqref{cls-mn}, we will have
\begin{align*}
	\dot{\tilde{V}}_2=&  \mathcal{K} \Big[ 2 \eta \tilde{x}^T  \tilde{v}+\tilde{v}^T \tilde{v} +\tilde{x}^T f(x,v,t) - \eta \tilde{x}^T   (\tilde{x}+\tilde{v}) \\
	&-\eta \tilde{x}^T (\psi' \otimes I_{\mathfrak{p}}) \mbox{sgn} (\tilde{x}+\tilde{v})  +\tilde{v}^T f(x,v,t)\\
	& - \eta \tilde{v}^T  (\tilde{x}+\tilde{v}) - \eta \tilde{v}^T ( \psi' \otimes I_{\mathfrak{p}}) \mbox{sgn} (\tilde{x}+\tilde{v})  \Big]   \\
	=&   \Big \{-\eta \tilde{x}^T \tilde{x} -(\eta-1)\tilde{v}^T \tilde{v} +(\tilde{x}^T + \tilde{v}^T) f(x,v,t)  \\
			&-\eta (\tilde{x}^T + \tilde{v}^T) (\psi' \otimes I_{\mathfrak{p}}) \mbox{sgn} (\tilde{x}+\tilde{v})  \Big \},
\end{align*}
where we have used the fact that $\mathcal{K} \Big[(\tilde{x}^T + \tilde{v}^T) (\psi' \otimes I_{\mathfrak{p}}) \mbox{sgn} (\tilde{x}+\tilde{v}) \Big ] =\mathcal{K} \Big[\sum\limits_{i=1}^{n}  \psi'_i (\tilde{x}_i+\tilde{v}_i)^T\mbox{sgn}(\tilde{x}_i+\tilde{v}_i) \Big ]= \mathcal{K} \Big[ \sum\limits_{i=1}^{n} \psi_i' \| \tilde{x}_i + \tilde{v}_i \|_1 \Big ]=\Big \{ \sum\limits_{i=1}^{n} \psi_i' \| \tilde{x}_i + \tilde{v}_i \|_1 \Big \}$.
Note that the set-valued Lie derivative of $V_2$ is a
singleton and the function $V_2$ is continuously differentiable.
It follows from $\dot{V}_2 \in \dot{\tilde{V}}_2$, where $\dot{V}_2$
denotes the derivative of $V_2$ , that
\begin{align*}
	\dot{V}_2	=&  -\eta \tilde{x}^T  \tilde{x}-(\eta-1) \tilde{v}^T  \tilde{v} +\sum\limits_{i=1}^{n} (\tilde{x}_i+\tilde{v}_i)^T f_i(x_i,v_i,t)  \\
	&- \eta \sum\limits_{i=1}^{n} \psi'_i \big \|\tilde{x}_i+\tilde{v}_i \big \|_1\\
	=&  -\eta \tilde{x}^T  \tilde{x}-(\eta-1) \tilde{v}^T  \tilde{v}  \\
	&+\sum\limits_{i=1}^{n} (\tilde{x}_i+\tilde{v}_i)^T (f_i(x_i,v_i,t)-f_i(\mathbf{0}_{\mathfrak{p}},\mathbf{0}_{\mathfrak{p}},t)) \\
	&+\sum\limits_{i=1}^{n} (\tilde{x}_i+\tilde{v}_i)^T f_i(\mathbf{0}_{\mathfrak{p}},\mathbf{0}_{\mathfrak{p}},t) \\
	&- \eta \sum\limits_{i=1}^{n} (\|x_i\|_1+\|v_i\|_1+\gamma) \big \|\tilde{x}_i+\tilde{v}_i \big \|_1 \\
		\leq &  -\eta \tilde{x}^T  \tilde{x}-(\eta-1) \tilde{v}^T  \tilde{v}  \\ 
		&- \eta \sum\limits_{i=1}^{n} (\|x_i\|_1+\|v_i\|_1+\gamma) \big \|\tilde{x}_i+\tilde{v}_i \big \|_1 \\
	&+\sum\limits_{i=1}^{n}  \Big (\rho_1 \|x_i\|_1+\rho_2 \|v_i\|_1 +\rho_3 +\rho_4   \Big) \big \|\tilde{x}_i+\tilde{v}_i \big \|_1 \\
	\leq &  -\eta \tilde{x}^T  \tilde{x}-(\eta-1) \tilde{v}^T  \tilde{v} <0,
\end{align*}
where we have used Assumption \ref{lip} to obtain the first inequality and $\eta >\max \{ 1, \rho_1,\rho_2 \}$ and $\gamma > \rho_3+\rho_4$ to obtain the second inequality.
Therefore, by using Theorem 4.10 in \cite{khalil2002nonlinear}, it is concluded that $\begin{bmatrix}
\tilde{x} \\
\tilde{v}
\end{bmatrix}=\mathbf{0}_{2n\mathfrak{p}}$ is globally exponentially stable. 
Thus, using \eqref{pi-avrri},  $\lim\limits_{t \to \infty} x_i = \frac{1}{n} \sum\limits_{j=1}^{n} r_j$ and 
$ \lim\limits_{t \to \infty} v_i = \frac{1}{n} \sum\limits_{j=1}^{n} v^r_j$.
\endproof
\begin{remark}
As it can be seen, by using algorithm \eqref{p-term-adap}-\eqref{con-inp-adap}, each agent can achieve the distributed average tracking, where there is no constraint on the input signals and the nonlinear terms in the agents' dynamics are unknown and heterogeneous.
Due to the presence of the unknown term $f_i(\cdot,\cdot,t)$ in the agents' dynamics, the existing algorithms for double-integrator agents are not applicable to achieve the distributed average tracking. 
For example, by employing the algorithm in \cite{fei2015DAC} for \eqref{non-agents-dynamic}, the two equalities $\sum_{j=1}^n x_j =\sum_{j=1}^n r_j$ and $\sum_{j=1}^n v_j = \sum_{j=1}^n v_j^r$ do not hold anymore.
In fact, the unknown term $f_i(\cdot,\cdot,t)$ functions as a disturbance and will not allow the average of the positions and velocities to track the average of the input signals and the average of the input velocities, respectively.
This shows the essence of using the local filter \eqref{p-term-adap} in our algorithm.

\end{remark}
\begin{remark} 
	In the proposed algorithm \eqref{p-term-adap}-\eqref{con-inp-adap}, correct position and velocity initialization are not required, which makes the algorithm applicable for platforms with physical agents. Note that  in real applications the correct initialization for physical variables might no be feasible.
\end{remark}
\subsection{Distributed Average Tracking for Physical Double-Integrator Agents}\label{section4}
In the proposed algorithm in Subsection \ref{section2}, the agents are described by heterogeneous nonlinear second-order dynamics, where the nonlinear term $f_i(\cdot,\cdot,t)$ satisfies a Lipschitz-like condition. 
However, in some applications the agents' dynamics can be linearized as double-integrator dynamics.
Therefore, in this subsection we modify the proposed algorithm in Subsection \ref{section2} for a group of agents  with double-integrator dynamics,
\begin{align} \label{agents-dynamic}
\dot{x}_i=v_i, \qquad \dot{v}_i= u_i,  \qquad  i=1,\ldots,n,
\end{align}
where $x_i$, $v_i$ and $u_i$ are introduced in Subsection \ref{section2}.
Since the nonlinear term $f_i(\cdot,\cdot,t)$ does not exist in the agents' dynamics, the local filter is not required here anymore.
Therefore, we can directly design $u_i$ to drive the agents' positions and velocities to track the average of  input signals and input velocities, respectively. 
	The control input for agent $i$, $i=1,\ldots,n$, is designed as
	\begin{align}\label{cont-lin} 
	u_i=&-\kappa (x_i-x^r_i) -\kappa (v_i-v^r_i) \notag \\
	&-\beta \sum\limits_{j=1}^{n} a_{ij}    (\psi_{i}+\psi_{j}) \mbox{sgn} \big[
	(x_i+v_i)-(x_j+v_j) \big] +a^r_i, 
	\end{align}
	where $\psi_{i}$ is defined in Subsection \ref{section2} and $\beta$, $\kappa$, $\gamma \in \mathbb{R}^+$ will be designed later.
	\begin{theorem} \label{DAC-lin}
		Under the control input given by  \eqref{cont-lin} for system \eqref{agents-dynamic}, the distributed average tracking goal \eqref{goal} is achieved asymptotically, provided that Assumption \ref{conn-graph} holds, $\beta>\kappa$ and $\kappa >\max \{1, \frac{\lambda^2_{\max} (L)}{2 \lambda_2^2(L)}\}$.
	\end{theorem}
	\emph{Proof}:
		Here the proof is very similar to the first step of Theorem \ref{DAC-lip-adap-adap} proof and hence the detail is omitted here.
		First, the following Lyapunove function is employed to prove that $\lim\limits_{ t \to \infty} x_i = \frac{1}{n} \sum\limits_{j=1}^{n} x_j$ and $\lim\limits_{ t \to \infty} v_i = \frac{1}{n} \sum\limits_{j=1}^{n} v_j$,  $i=1,\ldots,n$,
			\begin{align*}
				V=& \frac{1}{2} 
				\begin{bmatrix}
					e_x^T && e_v^T
				\end{bmatrix}
				\big ( L \otimes
				\begin{bmatrix}
					2\kappa  I_{\mathfrak{p}}  &&  I_{\mathfrak{p}}  \\
					I_{\mathfrak{p}}   && I_{\mathfrak{p}}
				\end{bmatrix}
				\big )
				\begin{bmatrix}
					e_x \\
					e_v
				\end{bmatrix},
			\end{align*}
		where $e_x=(M \otimes I_\mathfrak{p})x$ and $e_v=(M \otimes I_\mathfrak{p})v$ and $x$, $v$, $M$ are defined in Subsection \ref{section2}.
%
Second, it is shown that $\lim\limits_{t \to \infty} \sum\limits_{j=1}^{n} x_j=  \sum\limits_{j=1}^{n} r_j$ and $ \lim\limits_{t \to \infty} \sum\limits_{j=1}^{n} v_j=\sum\limits_{j=1}^{n} v^r_j$ asymptotically.
To prove that, we employ the control input \eqref{cont-lin} for \eqref{agents-dynamic} and rewrite the closed-loop dynamics as
	\begin{align}\label{sum-cl-lin}
		\dot{S}_x=& S_v, \notag \\
		\dot{S}_v=&- \kappa S_x- \kappa S_v,
	\end{align}
	where ${S_x=\sum_{j=1}^n (x_j- r_j)}$ and ${S_v=\sum_{j=1}^n (v_j-v_j^r) }$.
By using the same analysis as Theorem \ref{DAC-lip-adap-adap} and employing the results of the first step, it is finally concluded that $\lim\limits_{t \to \infty} x_i = \frac{1}{n} \sum_{j=1}^n r_j$ and ${\lim\limits_{t \to \infty} v_i = \frac{1}{n} \sum_{j=1}^n v_j^r}$.
	\endproof
\begin{remark}		
The introduced algorithm in \cite{fei2015DAC} can achieve the distributed average tracking provided that $a^r_i$ is bounded. 
		However, the algorithm \eqref{cont-lin} is more general and solves the problem regardless of any constraint on $a^r_i$. 
	\end{remark}
	\begin{corollary}
		Suppose that each agent is described by the following heterogeneous double-integrator dynamics
		\begin{align*}
			\dot{x}_i=& v_i, \notag \\
			m_i \dot{v}_i=& u_i, \qquad 0 < m_i \leq \bar{m}_i, \qquad i=1,\ldots,n,
		\end{align*}
		where $\bar{m}_i \in \mathbb{R}^+$.
		Let $\bar{m}=\max\limits_{i=1,\ldots,n} \bar{m}_i$.
		If the input acceleration $a^r_i$ in \eqref{ref-dynamic} is arbitrary, by employing the local filter and the control input defined in, respectively, \eqref{p-term-adap} and \eqref{con-inp-adap}, where $\beta>\kappa$, $\kappa >\max \{1, \frac{\lambda^2_{\max} (L)}{2 \lambda_2^2(L)}\}$, $\gamma > 0$, and $\eta> \max \{ 1, \frac{\bar{m}}{2} \}$, the distributed average tracking is achieved.
%
		
		\emph{Proof}:
		The proof is similar to the proof of Theorem \ref{DAC-lip-adap-adap}, where \eqref{LypV2} is replaced with
		\begin{align*}\label{LypV3}
			V_2= \frac{1}{2} \sum\limits_{i=1}^{n}
			\begin{bmatrix}
				\tilde{x}_i^T && \tilde{v}_i^T
			\end{bmatrix} 
			\begin{bmatrix}
				\frac{2\eta}{m_i} I_{n\mathfrak{p}} && I_{n\mathfrak{p}} \\
				I_{n\mathfrak{p}} && I_{n\mathfrak{p}}
			\end{bmatrix}
			\begin{bmatrix}
				\tilde{x}_i \\
				\tilde{v}_i
			\end{bmatrix}.
		\end{align*}
				\endproof
	\end{corollary}

	\section{Simulation}
\begin{figure}
		\centering
		\includegraphics[width=0.5\textwidth]{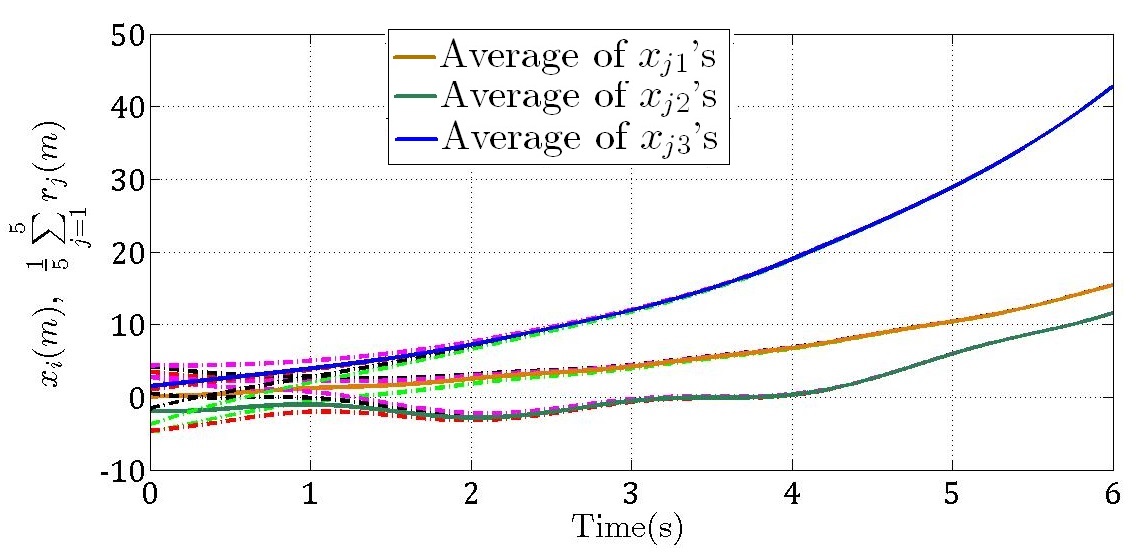}\vspace{-0.15cm}
		\caption{The positions of $4$ nonlinear agents, where $f_i(\cdot,\cdot,t)$ satisfy Assumption \ref{lip}.
			The solid lines and the
			dashed lines describe, respectively the average of  input signals and the position of agents.}
		\label{xb}
	\end{figure}

\begin{figure}
		\centering
		\includegraphics[width=0.5\textwidth]{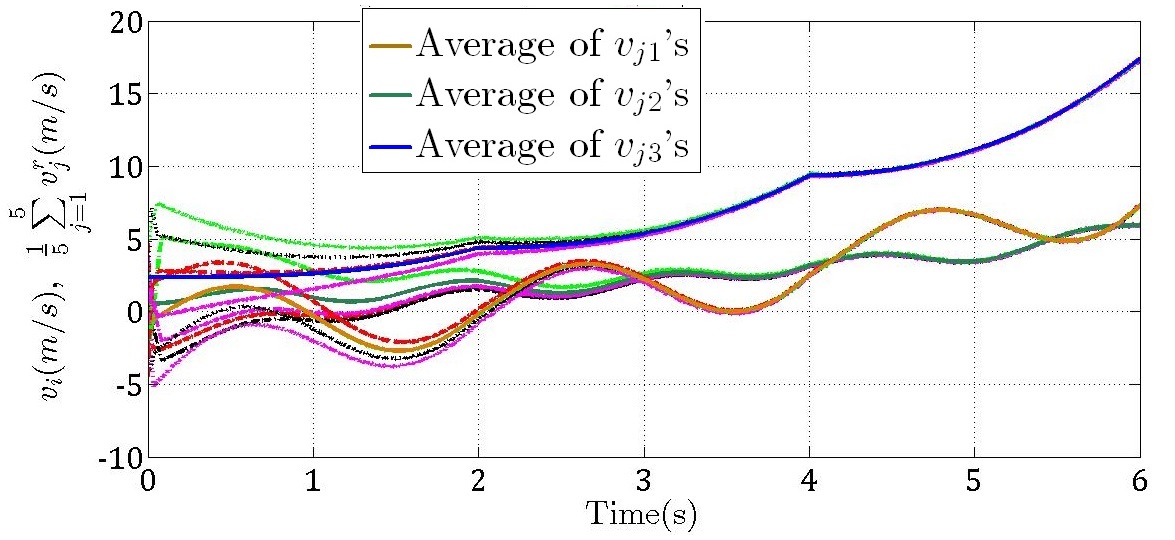}\vspace{-0.15cm}
		\caption{The velocities of $4$ agents and the average of  input velocities.
			The solid lines and the
			dashed lines describe, respectively the average of input velocities and the velocity of the agents.}
		\label{vb}
	\end{figure}

	Numerical simulation results are given in this section to illustrate the effectiveness of the theoretical results obtained in Subsection \ref{section2}.
	It is assumed that there are four agents $(n = 5)$.
	The nonlinear term $f_i(\cdot,\cdot,t)$ for agent $i$ is chosen as [26]
	\begin{align*}
		f_i(x,y,t)=i \times
		\begin{bmatrix}
			0 \\
			0\\
			- \beta \epsilon \sin (\omega x_{1})
		\end{bmatrix}
		+
		\begin{bmatrix}
			\delta ( y_{2} -y_{1} h(y_{1})) \\
			y_{1} - y_{2} +y_{3} \\
			-\beta y_{2} - \mu y_{3}
		\end{bmatrix},
	\end{align*}
	where $\delta =10$, $\beta = 19.53$, $\mu=0.1636$, $\epsilon=0.2$, $\omega = 0.5$ and $h(y_{1})=-0.7831 y_{1}- 0.324 (|y_{1} +1|-| y_{1}-1|)$. 
	It is easy to verify that the above nonlinear functions satisfy
	Assumption \ref{lip}.
	The input acceleration for agent $i$, $i = 1, \ldots,4$, is given by
	$a^r_i(t)=i \times [\sin(5t)+0.1 t \times \mbox{mod}(t,2),3\cos(3t)+0.2 t \times \mbox{mod}(t,2),0.3 t \times \mbox{mod}(t,2)]^T$.
	The initial positions and velocities of the agents are set randomly within the range $[-10,10]$.
	We denote the $j$th component of $x_i$ as $x_{ij}$. Similar notations are used for $v_i$, $x^r_i$, and $v_i^r$.
	The control parameters for all agent are chosen as $\beta=3$, $\kappa=1.5$, $\eta=5$ and $\gamma=0.001$.
	We simulate the algorithm defined by \eqref{p-term-adap}-\eqref{con-inp-adap}. Fig. \ref{xb} shows the positions of the agents and the average of the  input signals. Clearly, all agents have tracked the average of the  input signals in the presence of the nonlinear term $f_i(\cdot,\cdot,t)$ in the agents' dynamics. Fig. \ref{vb} shows the velocities of the agents and the average of the input velocities. We see that the distributed average tracking is achieved for the agents' velocities too.			
	
	\section{CONCLUSIONS}
	In this paper, distributed average tracking of physical second-order agents with heterogeneous nonlinear dynamics was studied, where there is no constraint on input signals. 
	The nonlinear terms in agents' dynamics satisfy the a Lipschitz-like condition that is more  general than the Lipschitz condition.
	For each agent, a control input combined with a local filter was designed. 
	The idea is that first the filter's outputs converge to the average of the input signals and input velocities asymptotically and then the agent's position and velocity are driven to track the filter outputs. 
	Since the nonlinear terms can be unbounded, a state-dependent time-varying gain was introduced in each filter's dynamics.
	Then, the algorithm was modified to achieve the distributed average tracking for physical second-order agents.
	In this algorithm, the filter is not required and a novel state-dependent time-varying gain was designed to solve the problem when there is no constraint on input signals.

	\bibliographystyle{IEEEtran}
	\bibliography{refs}

\end{document}